\documentclass[12pt]{amsart}

\advance\hoffset-15truemm

\let\Cal\mathcal
\let\Bbb\mathbb
\let\frak\mathfrak
\let\phi\varphi

\newcommand{\x}{\times}
\renewcommand{\o}{\circ}
\newcommand{\al}{\alpha}

\newcommand{\Ga}{\Gamma}
\newcommand{\ga}{\gamma}
\newcommand{\La}{\Lambda}
\newcommand{\om}{\omega}
\newcommand{\Om}{\Omega}
\newcommand{\la}{\lambda}
\newcommand{\ph}{\phi}
\newcommand{\Ph}{\Phi}
\newcommand{\ps}{\psi}
\newcommand{\ka}{\kappa}
\renewcommand{\th}{\theta}
\newcommand{\ze}{\zeta}

\newcommand{\tca}{{\tilde{\Cal A}}}

\newcommand{\tcs}{{\tilde{\Cal S}}}
\newcommand{\tct}{{\tilde{\Cal T}}}
\newcommand{\tcg}{{\tilde{\Cal G}}}
\newcommand{\Aut}{\operatorname{Aut}}
\newcommand{\aut}{\operatorname{aut}}
\newcommand{\Ad}{\operatorname{Ad}}

\newcommand{\gr}{\operatorname{gr}}
\newcommand{\im}{\operatorname{im}}

\newcommand{\fg}{{\mathfrak g}}
\newcommand{\fp}{{\mathfrak p}}
\newcommand{\tg}{{\tilde{\mathfrak g}}}
\newcommand{\tp}{{\tilde{\mathfrak p}}}

\newcommand{\xxb}{%
\begin{picture}(38,6)\put(3,3){\line(1,0){28}}%
\put(3,3){\makebox(0,0){$\x$}}\put(18,3){\makebox(0,0){$\x$}}%
\put(33,3){\makebox(0,0){$\o$}}\end{picture}}

\newcommand\xbdbb{\begin{picture}(76,6)\put(3,3){\line(1,0){13}}%
\put(20,3){\line(1,0){6}}\put(61,3){\line(-1,0){11}}%
\put(46,3){\line(-1,0){6}}\put(34,3){\makebox(0,0){\dots}}%
\put(3,3){\makebox(0,0){$\x$}}\put(18,3){\makebox(0,0){$\o$}}%
\put(48,3){\makebox(0,0){$\o$}}\put(63,3){\makebox(0,0){$\o$}}\end{picture}}

\newcommand\Cxbdbb{\begin{picture}(76,6)\put(3,3){\line(1,0){13}}%
\put(20,3){\line(1,0){6}}\put(61,2){\line(-1,0){11}}%
\put(61,4){\line(-1,0){11}}%
\put(46,3){\line(-1,0){6}}\put(34,3){\makebox(0,0){\dots}}%
\put(3,3){\makebox(0,0){$\x$}}\put(18,3){\makebox(0,0){$\o$}}%
\put(48,3){\makebox(0,0){$\o$}}\put(63,3){\makebox(0,0){$\o$}}%
\put(57,3){\makebox(0,0){$<$}}\end{picture}}

\newtheorem*{prop*}{Proposition}

\newtheorem*{thm*}{Theorem}

\newtheorem*{lem*}{Lemma}

\newtheorem*{kor*}{Corollary}

\theoremstyle{definition}

\newtheorem*{definition*}{Definition}

\newtheorem*{example*}{Example}

\begin{document}

\title{Two constructions with parabolic geometries}
\author{Andreas \v Cap}
\thanks{supported by project P15747-N05 of the ``Fonds zur
F\"orderung der wissenschaftlichen Forschung'' (FWF)}
\date{April 19, 2005}

\address{Institut f\"ur Mathematik, Universit\"at Wien,
Nordbergstra\ss e~15, A--1090 Wien, Austria, and International Erwin
Schr\"odinger Institute for Mathematical Physics, Boltzmanngasse 9,
A-1090 Wien, Austria}

\email{andreas.cap@esi.ac.at}

\begin{abstract}
This is an expanded version of a series of lectures delivered at the
25th Winter School ``Geometry and Physics'' in Srni. 

After a short introduction to Cartan geometries and parabolic
geometries, we give a detailed description of the equivalence between
parabolic geometries and underlying geometric structures. 

The second part of the paper is devoted to constructions which relate
parabolic geometries of different type. First we discuss the
construction of correspondence spaces and twistor spaces, which is
related to nested parabolic subgroups in the same semisimple Lie
group. An example related to twistor theory for Grassmannian
structures and the geometry of second order ODE's is discussed in
detail. 

In the last part, we discuss analogs of the Fefferman construction,
which relate geometries corresponding different semisimple Lie groups.
\end{abstract}


\maketitle

\section{Introduction}\label{1}
This is an expanded version of a series of plenary lectures at the
25th Winter School ``Geometry and Physics'' in Srni. I would like to thank
the organizers for giving me the opportunity to present this series.

The concept which is nowadays known as a Cartan geometry was
introduced by E.~Cartan under the name ``generalized space'' in order
to build a bridge between geometry in the sense of F.~Klein's
Erlangen program and differential geometry. This concept associates to
an arbitrary homogeneous space $G/H$ the notion of a \textit{Cartan
  geometry of type $(G,H)$}, which is a differential geometric
structure on smooth manifolds whose dimension equals the dimension of
$G/H$. A manifold endowed with such a geometry can be considered as a
``curved analog'' of the homogeneous spaces $G/H$. Although Cartan
geometries are an extremely general concept, there are several
remarkable results which hold for all of them, see \ref{2.2}.

The most interesting examples of Cartan geometries are those, in which
the Cartan geometry is equivalent to some simpler underlying
structure. Obtaining the Cartan geometry from the underlying structure
usually is a highly nontrivial process which often involves
prolongation. Cartan himself found many examples of this situation,
ranging from conformal and projective structures via 3--dimensional CR
structures to generic rank two distributions in manifolds of dimension
five.

Parabolic geometries are Cartan geometries of type $(G,P)$, where $G$
is a semisimple Lie group and $P\subset G$ is a parabolic
subgroup. The corresponding homogeneous spaces $G/P$ are the
so--called generalized flag manifolds which are among the most
important examples of homogeneous spaces. Under the conditions of
regularity and normality, parabolic geometries always are equivalent to
underlying structures. This basically goes back to the pioneering
works of N.~Tanaka, see e.g.~\cite{Tanaka79}.

In section \ref{2} of this article we give a precise description of
the underlying structures which are equivalent to regular normal
parabolic geometries. In this underlying picture, the structures are
very diverse, including in particular the four examples of structures
listed above. From that point of view, parabolic geometries offer a
unified approach to a broad variety of geometric structures.

Some of the advantages of this unified approach will be discussed in
the remaining two sections. They are devoted to constructions which
relate parabolic geometries of different types. The common feature of
these constructions is that they are quite transparent in the picture
of Cartan geometries, while from the point of view of the underlying
structures they are often surprising. 

Section \ref{3} is devoted to the construction of correspondence
spaces, which is associated to nested parabolic subgroups in one
semisimple Lie group. Trying to characterize the geometries obtained
in that way, one is lead to the notion of a twistor space and obtains
several classical examples of twistor theory. In the end one arrives
at a complete local characterization of correspondence spaces in terms
of the harmonic curvature. We give a detailed discussion of one
example of this situation related to the geometry of systems of second
order ODE's.

The last section is devoted to Fefferman's construction of a conformal
structure on the total space of a circle bundle over a CR manifold and
analogs of this construction. From the point of view of Cartan
geometries, the basic input for these constructions is an inclusion
$i:G\to\tilde G$ between semisimple groups together with appropriately
chosen parabolic subgroups $P\subset G$ and $\tilde P\subset\tilde
G$. Then the construction relates geometries of type $(G,P)$ to
geometries of type $(\tilde G,\tilde P)$.

\section{Cartan geometries and parabolic geometries}\label{2}
We start with some general background on Cartan geometries. 

\subsection{Homogeneous spaces and the Maurer Cartan form}\label{2.1}
Let $G$ be any Lie group and let $H\subset G$ be a closed
subgroup. The basic idea behind Cartan geometries is to endow the
homogeneous space $G/H$ with a geometric structure, whose automorphisms
are exactly the left actions of the elements of $G$. The natural
projection $G\to G/H$ is well known to be a principal bundle with
structure group $H$. Left multiplication by $g\in G$ lifts the action
of $g$ on $G/H$ to an automorphism of this principal bundle. Of
course, the group of principal bundle automorphisms of $G\to G/H$ is
much bigger than just the left translations, so an additional
ingredient is needed to recognize left translations. 

It turns out that the right ingredient is the (left) \textit{Maurer
Cartan form} $\om^{MC}\in\Om^1(G,\frak g)$. Recall that this is just a
different way to encode the trivialization of the tangent bundle $TG$
by left translations. By definition, for $\xi\in T_gG$ we have
$$
\om^{MC}(\xi)=T\la_{g^{-1}}\cdot\xi\in T_eG=\frak g,
$$
where $\la_{g^{-1}}$ denotes left translation by $g^{-1}$.

\begin{prop*}
Let $G$ be a Lie group and let $H\subset G$ be a closed subgroup such
that the homogeneous space $G/H$ is connected. Then the left
translations $\la_g$ are exactly the principal bundle automorphisms of
$G\to G/H$ which pull back $\om^{MC}$ to itself. 
\end{prop*}

For later use, we note some further properties of $\om^{MC}$. As we
have noted above, $\om(L_X)=X$ for all $X\in\frak g$. Note that for
$X\in\frak h$, the vector field $L_X$ coincides with the fundamental
vector field $\ze_X$ on the principal bundle $G\to G/H$ generated by
$X$. For $g\in H$, consider the right translation $r^g$ by $g$. Using
that the adjoint action of $g$ is the derivative of the conjugation by
$g$, one immediately verifies that
$(r^g)^*\om^{MC}=\Ad(g^{-1})\o\om^{MC}$. Note that for $g\in H$, the
map $r^g$ is the principal right action on the bundle $G\to G/H$.
Finally, there is the Maurer--Cartan equation: The fact that
$[L_X,L_X]=L_{[X,Y]}$ for all $X,Y\in\frak g$ implies that
$d\om^{MC}(\xi,\eta)+[\om(\xi),\om(\eta)]=0$ for all vector fields
$\xi$ and $\eta$ on $G$.

\subsection{Cartan geometries}\label{2.2}
The definition of a Cartan geometry is now obtained by replacing $G\to
G/H$ by an arbitrary principal $H$--bundle and $\om^{MC}$ by a form
which has all the properties of $\om^{MC}$ that make sense in the more
general setting.

\begin{definition*}
(1) A \textit{Cartan geometry} of type $(G,H)$ on a smooth manifold
$M$ is a principal $H$--bundle $p:\Cal G\to M$ together with a one
form $\om\in\Om^1(\Cal G,\frak g)$ (the \textit{Cartan connection})
such that
\begin{center}
\parbox{0.9\textwidth}{
$\bullet$ $(r^h)^*\om=\Ad(h)^{-1}\o\om$ for all $h\in H$.
\newline
$\bullet$ $\om(\ze_A)=A$ for all $A\in\frak h$.
\newline
$\bullet$ $\om(u):T_u\Cal G\to\frak g$ is a linear isomorphism for all
$u\in\Cal G$.}
\end{center}
(2) A \textit{morphism} between two Cartan geometries $(\Cal G\to
M,\om)$ and $(\tcg\to\tilde M,\tilde\om)$ is a principal bundle
homomorphism $\Ph:\Cal G\to\tcg$ such that $\Ph^*\tilde\om=\om$.

\noindent
(3) The \textit{curvature} $K\in\Om^2(\Cal G,\frak g)$ of a Cartan
    geometry $(\Cal G\to M,\om)$ of type $(G,H)$ is defined by 
$$
K(\xi,\eta)=d\om(\xi,\eta)+[\om(\xi),\om(\eta)],
$$
for $\xi,\eta\in\Cal X(\Cal G)$. 
\end{definition*}

\smallskip

Notice that a Cartan geometry is a local structure, i.e.~it can be
restricted to open subsets: For $(p:\Cal G\to M,\om)$ and an open
subset $U\subset M$, we simply have the restriction $(p:p^{-1}(U)\to
U,\om|_{p^{-1}(U)})$. The curvature evidently is a local invariant,
i.e.~the curvature of this restricted geometry is the restriction of
the original curvature.

Any morphism $\Ph$ between two Cartan geometries as in (2) has an
underlying smooth map $\ph:M\to\tilde M$. It turns out (see
\cite[chapter 5]{Sharpe}) that $\Ph$ is determined by $\ph$ up to a smooth
function from $M$ to the maximal normal subgroup of $G$ which is
contained in $H$. In all cases of interest, this subgroup is trivial
or at least discrete, whence this map has to be locally constant. In
fact, it is necessary to included the possibility of having various
morphisms covering the same base map to deal with structures analogous
to Spin structures.

By definition $(G\to G/H,\om^{MC})$ is a Cartan geometry of type
$(G,H)$, and Proposition \ref{2.1} exactly tells us that the
automorphisms of this geometry are exactly the left translations by
elements of $G$. This geometry is called the \textit{homogeneous
model} of Cartan geometries of type $(G,H)$. 

The Maurer--Cartan equation noted in the end of \ref{2.1} exactly says
that the curvature of the homogeneous model vanishes identically.
Indeed, the curvature exactly measures to what extent the
Maurer--Cartan equation fails to hold. One of the nice features of
Cartan geometries is that vanishing of the curvature characterizes the
homogeneous model locally, i.e.~any Cartan geometry of type $(G,H)$
with vanishing curvature is locally isomorphic to $(G\to
G/H,\om^{MC})$, see \cite[chapter 5]{Sharpe}. More generally, the
curvature (at least in principle) provides a solution to the
equivalence problem. This is one of the reasons why already
associating to some geometric structure a canonical Cartan connection
is a powerful result. For the main part of the theory of parabolic
geometries however, the existence of a canonical Cartan connection is
only the starting point.

There are other interesting features of general Cartan
geometries, for example:
\begin{itemize}
\item For any Cartan geometry $(p:\Cal G\to M,\om)$ of type $(G,H)$,
the automorphism group $\Aut(\Cal G,\om)$ is a Lie group of dimension
$\leq\dim(G)$. The Lie algebra $\aut(\Cal G,\om)$ can be described
completely, and analyzing its algebraic structure leads to interesting
results, see \cite{Srni04}. 
\item The homogeneous model $(G\to G/H,\om^{MC})$ satisfies a Liouville
type theorem. If $U$ and $V$ are open subsets of $G/H$ then any
isomorphism between the restrictions of the geometry to these open
subsets uniquely extends to an automorphism of the homogeneous model. 
\item There are various general tools available for Cartan geometries,
for example the notions of distinguished curves and of normal
coordinates. 
\end{itemize}

\subsection{Cartan geometries determined by underlying
  structures}\label{2.3} 
The results listed above become particularly powerful if a Cartan
geometry is obtained as an equivalent description of some underlying
geometric structure. A very simple example is provided by Riemannian
geometries, which correspond to the case that $G$ is the Euclidean
group $\text{Euc}(n)$ and $H=O(n)$. The Lie algebra $\frak g$ is
isomorphic to $\frak h\oplus\Bbb R^n$ as an $H$--module. Therefore, a
Cartan connection of type $(G,H)$ on a principal $O(n)$--bundle $\Cal
G\to M$ decomposes into an $\Bbb R^n$--valued form $\th$ and a $\frak
h$--valued form $\ga$ which both are $H$--equivariant. Then $\th$
defines a reduction of the linear frame bundle of $M$ to the structure
group $O(n)$, which is equivalent to a Riemannian metric on $M$. The
form $\ga$ defines a principal connection on $\Cal G$ which is
equivalent to a metric connection on $M$. If $\ga$ is torsion free,
then it must be the Levi--Civita connection. Conversely, starting from
a Riemannian manifold, one obtains a torsion free Cartan geometry by
using the orthonormal frame bundle endowed with the soldering form and
the Levi--Civita connection. In that way, one obtains an equivalence
between torsion free Cartan geometries of type $(G,H)$ and
$n$--dimensional Riemannian manifolds.

The results discussed above then imply
\begin{itemize}
\item The isometry group of any Riemannian manifold is a Lie group of
dimension $\leq \tfrac{1}{2}\dim(M)(\dim(M)-1)$.
\item Any isometry between two open subsets of Euclidean space is the
restriction of a uniquely determined Euclidean motion. 
\item The concepts of geodesics and Riemann normal coordinates for
Riemannian manifolds. 
\end{itemize}

The case of Riemannian metrics is rather easy, since the bundle $\Cal
G$ can be directly obtained from the underlying structure. In other
cases, one also has to construct this principal bundle, a process
which is usually called \textit{prolongation}. This also leads to
additional features. Let us discuss this in the case of conformal
structures, which is a model case for parabolic geometries:  

A conformal structure on a smooth manifold $M$ is given by an
equivalence class $[g]$ of Riemannian metrics on $M$. Here two metrics
$g$ and $\hat g$ are considered as equivalent if and only if $\hat
g=e^{2f} g$ for some smooth function $f$ on $M$. Equivalently, a
conformal structure can be defined as a reduction of structure group
of the frame bundle $\Cal PM$ to the group $CO(n)$ of conformal
isometries of $\Bbb R^n$.

It is a classical result of E.~Cartan, see \cite{Cartan:conf} that for
$n=\dim(M)\geq 3$ conformal structures admit a canonical normal Cartan
connection. Consider the semisimple Lie group $G:=SO(n+1,1)$. This
naturally acts on $\Bbb R^{n+2}$ and the action preserves the null
cone. Fix a nonzero null vector $v$ and let $P\subset G$ be the
stabilizer of the line $\Bbb Rv$. Then $P$ is an example of a
\textit{parabolic subgroup} of the semisimple Lie group $G$. It turns
out that $P$ contains an Abelian normal subgroup $P_+\cong\Bbb R^n$
such that $P/P_+=:G_0\cong CO(n)$. 

The relation of these groups to conformal geometry is the following:
The group $G$ acts transitively on the space of null lines in $\Bbb
R^{n+2}$, which is easily seen to be isomorphic to $S^n$. Since by
definition $P$ is the stabilizer of one null line we get $G/P\cong
S^n$ and this identifies $G$ with the group of conformal isometries of
$S^n$ and $P$ with the group of conformal isometries fixing a point
$x_0\in S^n$. It turns out that the projection from $P$ to $G_0\cong
CO(n)$ is given by passing from a conformal isometry fixing $x_0$ to
its tangent map in $x_0$, see \cite{mike:srni} for more details.

Now consider a manifold $M$ of dimension $n$ endowed with a conformal
structure $[g]$. The corresponding reduction of structure group is a
$G_0$--principal bundle $p_0:\Cal G_0\to M$ endowed with a canonical
differential form $\th$ called the \textit{soldering form}. Cartan's
result states that this bundle can be canonically extended to a
principal bundle $p:\Cal G\to M$ with structure group $P$ and the
soldering form $\th$ can be extended to a Cartan connection
$\om\in\Om^1(\Cal G,\frak g)$. If one requires the Cartan connection
$\om$ to satisfy a normalization condition, then it is uniquely
determined.

Conversely, given a principal $P$--bundle $p:\Cal G\to M$ endowed with
a Cartan connection $\om\in\Om^1(\Cal G,\frak g)$, one obtains a
$G_0$--principal bundle $\Cal G_0:=\Cal G/P_+\to M$ and $\om$ induces
a soldering form on that bundle, thus giving rise to a conformal
structure on $M$. In the end one obtains an equivalence of categories
between conformal structures and Cartan geometries of type $(G,P)$. 

The additional feature provided by this is that one obtains new
geometric objects. Viewing a conformal structure as a reduction to the
structure group $CO(n)$ of the linear frame bundle, one obtains a
natural vector bundle associated to each representation of
$CO(n)$. Since $CO(n)$ is a quotient of $P$, this gives rise to a
representation of $P$ and the resulting vector bundles can also be
viewed as associated to the Cartan bundle $\Cal G$. But the group $P$
admits more general representations than those coming from $G_0$, and
these give rise to new natural vector bundles and thus new geometric
objects. A particularly interesting case is to consider restrictions
to $P$ of representations of $G$. This leads to the so--called tractor
bundles, see \cite{BEG, tractors}. 

In a series of pioneering papers in the 1960's and 70's culminating in
\cite{Tanaka79}, N.~Tanaka showed that for all semisimple Lie groups
and parabolic subgroups normal Cartan geometries are determined by
underlying structures. These results have been put into the more
general context of filtered manifolds in the work of T.~Morimoto (see
e.g.~\cite{Morimoto}) and a new version of the result tailored to the
parabolic case was given in \cite{Cap-Schichl}. Otherwise put, these
results show that these underlying structures (which seemingly are
very diverse) admit canonical Cartan connections. Our next aim is to
give a uniform description of the underlying structures. 

\subsection{Generalized flag manifolds}\label{2.4}
We first collect some background on the homogeneous models of parabolic
geometries. We will use elementary definitions, which avoid structure
theory of Lie algebras.  

\begin{definition*}
Let $\frak g$ be a semisimple Lie algebra. A \textit{$|k|$--grading}
on $\frak g$ is a vector space decomposition 
$$
\frak g=\frak g_{-k}\oplus\dots\oplus\frak g_0\oplus\dots\oplus\frak g_k
$$
such that $[\frak g_i,\frak g_j]\subset \frak g_{i+j}$ and such that
the subalgebra $\frak g_-:=\frak g_{-k}\oplus\dots\oplus\frak g_{-1}$
is generated by $\frak g_{-1}$.
\end{definition*}

For given $\frak g$ there is a simple complete description of such
gradings (up to isomorphism) in terms of structure theory. For complex
$\frak g$, they are in bijective correspondence with sets of simple
roots of $\frak g$ and hence are conveniently denoted by Dynkin
diagrams with cros\-ses. For real $\frak g$ there is a similar
description in terms of the Satake diagram.

Let us make this more explicit for the case $\frak
g=\frak{sl}(n+1,\Bbb K)$ for $\Bbb K=\Bbb R$ or $\Bbb C$. Up to
isomorphism, each $|k|$--grading is determined by a block
decomposition of matrices: One decomposes $\Bbb R^{n+1}$ into $k+1$
blocks of sizes $i_0,\dots,i_k$. The $\frak g_0$ consists of all block
diagonal matrices, and for $i>0$, the component $\frak g_i$
(respectively $\frak g_{-i}$) consists of those matrices, which only
have nonzero entries in the $i$th blocks above (respectively below)
the main diagonal. The corresponding Dynkin diagram is obtained as
follows: Look at the matrix entries in the first diagonal above the
main diagonal. The block in which they are contained either lies in
$\frak g_0$ or in $\frak g_1$. Use a dot in the first and a cross
in the second case and connect each entry with a line to its (one or
two) neighbors. More explicitly, consider $\frak{sl}(4,\Bbb K)$ with
blocks of sizes $1$, $1$, and $2$. Then one obtains a $|2|$--grading
of the form
$$
\begin{pmatrix}\frak g_0 & \frak g_1 &\frak g_2&\frak g_2\\
\frak g_{-1} & \frak g_0 & \frak g_1 & \frak g_1\\
\frak g_{-2} & \frak g_{-1} &\frak g_0 &\frak g_0\\
\frak g_{-2} & \frak g_{-1} &\frak g_0 &\frak g_0\\ \end{pmatrix}
$$
and the corresponding Dynkin (respectively Satake) diagram with
crosses is $\xxb$.  

Putting $\fg^i:=\fg_i\oplus\dots\oplus\fg_k$ we obtain a filtration
$\fg=\fg^{-k}\supset\dots\supset\fg^k$ of $\fg$ such that
$[\fg^i,\fg^j]\subset\fg^{i+j}$. In particular, $\fp:=\fg^0$ is a
subalgebra of $\fg$ and $\fp_+:=\fg^1$ is a nilpotent ideal in $\frak
p$ such that $\frak p=\fg_0\oplus\frak p_+$ is a semidirect sum. The
subalgebras $\fp$ obtained in that way are exactly the
\textit{parabolic} subalgebras of $\fg$ used in representation theory.
In the complex case, a subalgebra of $\frak g$ is parabolic if and
only if it contains a maximal solvable subalgebra (i.e.~a Borel
subalgebra) of $\frak g$. In the real case, parabolic subalgebras are
defined via complexification. 

Given a (not necessarily connected) Lie group $G$ with Lie algebra
$\fg$, it turns out that the normalizer $P:=N_G(\frak p)$ of $\frak p$
in $G$ has Lie algebra $\frak p$. This is the \textit{parabolic
  subgroup} of $G$ associated to the parabolic subalgebra $\frak
p\subset\frak g$. It turns out that for $g\in P$, the adjoint action
$\Ad(g):\fg\to\fg$ not only preserves the filtration component $\fg^0$
but all filtration components $\fg^i$. Indeed, the whole filtration can
be reconstructed algebraically from $\fg^0=\frak p$. Further, one
defines a closed subgroup $G_0\subset P$ as the set of those $g\in P$,
whose adjoint action even preserves the grading of $\fg$. Then $G_0$
is reductive and has Lie algebra $\frak g_0$. One shows that $\exp$
defines a diffeomorphism from $\fp_+$ onto a closed subgroup
$P_+\subset P$ and $P$ is the semidirect product of $G_0$ and $P_+$.

A \textit{generalized flag variety} is a homogeneous space $G/P$ for a
semisimple Lie group $G$ and a parabolic subgroup $P\subset G$. These
homogeneous spaces are always compact and for complex $G$ they are the
only compact homogeneous spaces of $G$. In the complex case, $G/P$
carries a K\"ahler metric. Generalized flag manifolds are among the
most important examples of homogeneous spaces. They show up in many
areas of mathematics.

\textit{Parabolic geometries} are Cartan geometries of type $(G,P)$
for $G$ and $P$ as above. In \ref{2.3} we have seen that for an
appropriate choice of $G$ and $P$, such a structure (satisfying an
additional normalization condition) is equivalent to a conformal
Riemannian structure. Under the conditions of regularity and
normality, general parabolic geometries equivalent to a certain
underlying structure. We will next describe how this underlying
structure is obtained.

\subsection{The filtration of the tangent bundle}\label{2.5}
We first show how a parabolic geometry $(p:\Cal G\to M,\om)$ of type
$(G,P)$ gives rise to a filtration of the tangent bundle $TM$. Define
the \textit{adjoint tractor bundle} $\Cal AM$ of $M$ as $\Cal AM:=\Cal
G\x_P\fg$. (This is an important example of the concept of tractor
bundles discussed in \ref{2.3}.) Then we have the $P$--invariant
filtration $\{\fg^i\}$ of $\fg$, which gives rise to a filtration
$$
\Cal AM=\Cal A^{-k}M\supset\Cal A^{-k+1}M\supset\dots\supset\Cal A^kM
$$
of the adjoint tractor bundle by smooth subbundles. The Lie bracket
on $\fg$ induces a tensorial map $\{\ ,\ \}:\Cal AM\x\Cal AM\to\Cal
AM$. In particular, each fiber of $\Cal AM$ is a filtered Lie algebra
isomorphic to $\fg$.

The Cartan connection $\om$ leads to an identification $TM\cong\Cal
G\x_P\fg/\fp$, with the action coming from the adjoint action. The
Killing form of $\fg$ induces a duality between this $P$--module and
$\fp_+=\fg^1$, so $T^*M\cong \Cal G\x_P\fp_+=\Cal A^1M$. Hence $T^*M$
is a bundle of nilpotent filtered Lie algebras. On the tangent bundle,
there are similar but more subtle structures: From above, we see that
$TM\cong\Cal AM/\Cal A^0M$, and we obtain an induced filtration
$TM=T^{-k}M\supset\dots\supset T^{-1}M$ of the tangent bundle by
putting $T^iM:=\Cal A^iM/\Cal A^0M$. The associated graded bundle is
$$
\gr(TM)=\gr_{-k}(TM)\oplus\dots\oplus\gr_{-1}(TM),
$$
where $\gr_i(TM)=T^iM/T^{i+1}M$. By construction, this implies that
$\gr_i(TM)\cong\Cal G\x_P\fg^i/\fg^{i+1}$. By definition,
the subgroup $P_+\subset P$ acts trivially on this quotient. Hence the
$P$--action factorizes over $P/P_+\cong G_0$ and as a $G_0$--module we
have $\fg^i/\fg^{i+1}\cong \fg_i$.

On the level of principal bundles, we observe that the subgroup
$P_+\subset P$ acts freely on $\Cal G$. Hence the quotient $\Cal
G_0:=\Cal G/P_+$ is a principal bundle over $M$ with structure group
$P/P_+=G_0$. The Cartan connection $\om$ induces a bundle map from
$\Cal G_0$ to the frame bundle of $\gr(TM)$ which defines a reduction
of structure group. In particular, $\gr_i(TM)\cong\Cal
G_0\x_{G_0}\fg_i$, which is a refined version of the identification of
the representation spaces above. Putting the components together, we
see that $\gr(TM)\cong\Cal G_0\x_{G_0}\fg_-$. The Lie bracket on
$\fg_-$ is $G_0$--invariant and hence gives rise to a tensorial map
$\{\ ,\ \}$ on $\gr(TM)$. Hence for each $x\in M$, the space
$\gr(T_xM)$ is a nilpotent graded Lie algebra isomorphic to $\fg_-$.

\subsection{Filtered manifolds and their symbol algebras}\label{2.6}
A crucial observation for the sequel is that under a weak condition, a
filtration of the tangent bundle of a manifold gives rise to the
structure of a nilpotent graded Lie algebra on the associated graded
of each tangent space.

A \textit{filtered manifold} is a smooth manifold $M$ together with a
filtration $TM=T^{-k}M\supset\dots\supset T^{-1}M$ of the tangent
bundle by smooth subbundles, which is compatible with the Lie bracket
of vector fields, i.e.~such that for $\xi\in\Ga(T^iM)$ and
$\eta\in\Ga(T^jM)$ one always has $[\xi,\eta]\in\Ga(T^{i+j}M)$. 

Let $q:T^{i+j}M\to T^{i+j}M/T^{i+j+1}M=\gr_{i+j}(TM)$ be the natural
map, and consider the operator
$\Ga(T^iM)\x\Ga(T^jM)\to\Ga(\gr_{i+j}(TM))$ defined by
$(\xi,\eta)\mapsto q([\xi,\eta])$. Since the indices of the filtration
components are always negative, the bundles $T^iM$ and $T^jM$ are
contained in $T^{i+j+1}M$, which implies that this operator is
bilinear over smooth functions. Therefore, it is induced by a tensor
$T^iM\x T^jM\to\gr_{i+j}(TM)$. If $\xi\in\Ga(T^{i+1}M)$, then
$[\xi,\eta]\in\Ga(T^{i+j+1}M)$ so the result of this tensor depends
only on the classes of $\xi$ in $\gr_i(TM)$ and $\eta\in\gr_j(TM)$.
Taking together the various components, we obtain a tensor $\Cal
L:\gr(TM)\x\gr(TM)\to\gr(TM)$ which is called the \textit{Levi
  bracket}. By construction, this makes each of the spaces $\gr(T_xM)$
into a nilpotent graded Lie algebra, called the \textit{symbol
  algebra} of the filtered manifold at the point $x$. Consider a local
isomorphism between filtered manifolds, i.e.~a local diffeomorphism
$f$ such that each of the maps $T_xf$ is an isomorphism of filtered
vector spaces. Then each $T_xf$ induces and isomorphism between the
associated graded spaces to the tangent spaces, which is easily seen
to be an isomorphism of the symbol algebras.

Therefore, the symbol algebra should be considered as the first order
approximation of a filtered manifold in a point, similarly to the
tangent space at a point of an ordinary manifold. The usual tangent
space (viewed as an Abelian Lie algebra) is recovered in the case of
the trivial filtration $T^{-1}M=TM$.

A priory, the isomorphism class of the symbol algebra may change from
point to point, but the case that all symbol algebras are isomorphic
to a fixed nilpotent graded Lie algebra $\frak a$ is of particular
interest. In this case, there is a natural frame bundle for the vector
bundle $\gr(TM)$ with structure group the group $\Aut_{\gr}(\frak a)$
of all automorphisms of the Lie algebra $\frak a$, which in addition
preserve the grading. This is the replacement for the usual frame
bundle of a smooth manifold, which is again recovered in the special
case of the trivial filtration. 

\subsection{Regularity and normality}\label{2.7}
Let $(p:\Cal G\to M,\om)$ be a parabolic geometry of type
$(G,P)$. Then we have the curvature $K\in\Om^2(\Cal G,\frak g)$ of
$\om$ as introduced in \ref{2.2}. The defining properties of $K$
easily imply that it is horizontal and $P$--equivariant, so it defines a
two--form $\ka$ on $M$ with values in the bundle $\Cal G\x_P\fg=\Cal
AM$. Hence the Cartan curvature can be viewed as a two form on $M$
with values in the adjoint tractor bundle. 

The geometry  $(p:\Cal G\to M,\om)$ is called \textit{regular} if and
only if the curvature $\ka$ has the property that
$\ka(T^iM,T^jM)\subset \Cal A^{i+j+1}M$ for all $i,j<0$. Otherwise
put, regularity means that the curvature is concentrated in positive
homogeneities. 

Recall that a Cartan geometry of type $(G,P)$ is called
\textit{torsion free}, if its curvature $K\in\Om^2(\Cal G,\frak g)$
actually has values in $\frak p\subset\frak g$. In the parabolic case,
this can be nicely reformulated as $\ka$ lying in the subspace
$\Om^2(M,\Cal A^0M)\subset \Om^2(M,\Cal AM)$. From this description,
it is evident that torsion free parabolic geometries are automatically
regular, so regularity can be viewed as a condition avoiding
particularly bad types of torsion. Note that the condition is vacuous
for $|1|$--gradings. 

The geometric meaning of the regularity condition is easy to describe
(and also easy to prove): 
\begin{prop*}
Let $(p:\Cal G\to M,\om)$ be a parabolic geometry of type $(G,P)$, let
$\{T^iM\}$ be the induced filtration of the tangent bundle, and let
$\{\ ,\ \}$ be the tensorial Lie bracket on $\gr(TM)$ introduced in
\ref{2.5}. 

Then the geometry is regular if and only if the filtration $\{T^iM\}$
makes $M$ into a filtered manifold such that the natural bracket on
each symbol algebra coincides with $\{\ ,\ \}$. In particular, each
symbol algebra is isomorphic to $\frak g_-$. 
\end{prop*}

For regular geometries, the bundle $\Cal G_0\to M$ from \ref{2.5}
nicely ties into the concepts for filtered manifolds. The adjoint
action of $G_0$ on $\fg_-$ is by Lie algebra automorphisms which
preserve the grading (by definition of $G_0$), so it defines a
homomorphism $G_0\to\Aut_{\gr}(\frak g_-)$. This homomorphism turns
out to be infinitesimally injective provided that none of the simple
ideal of $\frak g$ is contained in $\fg_0$. This condition is very
harmless, since simple ideals contained in $\fg_0$ can be left out
without problems, so we will assume throughout that it is
satisfied. As we have noted in \ref{2.6}, the group $\Aut_{\gr}(\frak
g_-)$ is the natural structure group for the vector bundle $\gr(TM)$
since each symbol algebra is isomorphic to $\fg_-$. The bundle $\Cal
G_0$ can thus be interpreted as the filtered manifold version of a
first order $G_0$--structure. 

Now we have collected the two structures underlying a regular
parabolic geometry of type $(G,P)$ that we will need in the sequel:
\begin{itemize}
\item A filtration $\{T^iM\}$ of the tangent bundle such that each
symbol algebra is isomorphic to $\fg_-$. 
\item A reduction of structure group of the associated graded
$\gr(TM)$ to the structure group $G_0\subset\Aut_{\gr}(\frak g_-)$. 
\end{itemize}
Similarly to the soldering form used for classical first order
structures, this reduction of structure group can be expressed by
certain partially defined differential forms on the bundle $\Cal G_0$.
This leads to the description of underlying structures used in
\cite{Cap-Schichl}.  The collection of these two underlying structures
is called a \textit{regular infinitesimal flag structure}, see
\cite{Weyl}.

Fixing the underlying regular infinitesimal flag structure still
leaves a lot of freedom for the Cartan connection $\om$, so we need an
additional \textit{normalization condition}: Recall the the cotangent
bundle $T^*M$ can be naturally viewed as $\Cal G\x_P\fp_+=\Cal
A^1M$. Hence it naturally is a bundle of nilpotent Lie algebras with
the restriction of the algebraic bracket $\{\ ,\ \}$ of $\Cal AM$. Now
for $\ell>0$ we define a tensorial operator
$\partial^*:\La^{\ell}T^*M\otimes\Cal AM\to\La^{\ell-1}T^*M\otimes\Cal
AM$ by
\begin{align*}
\partial^*(\al_1\wedge&\dots\wedge\al_{\ell}\otimes s):=
\sum_{i=1}^\ell(-1)^i\al_1\wedge\dots\wedge\widehat{\al_i}
\wedge\dots\wedge\al_{\ell}\otimes\{\al_i,s\}+\\
&\sum_{i<j}(-1)^{i+j}\{\al_i,\al_j\}\wedge\al_1\wedge\dots\wedge
\widehat{\al_i}\wedge\dots\wedge\widehat{\al_j}\wedge\dots\wedge
\al_{\ell}\otimes s
\end{align*}
for $\al_r\in T^*M$ and $s\in\Cal AM$, where as usual the hats denote
omission. This is the differential in the standard complex computing
Lie algebra homology. In particular, $\partial^*\o\partial^*=0$, and
the quotients $\ker(\partial^*)/\im(\partial^*)$ are the pointwise Lie
algebra homologies of the Lie algebras $T^*_xM$ with coefficients in
the modules $\Cal A_xM$.

The homology groups $H_*(\frak p_+,\frak g)$ are naturally
$P$--modules and it is easy to see that $P_+$ acts trivially, so they
are obtained by trivially extending the action of $G_0$. Hence the
above bundles $\ker(\partial^*)/\im(\partial^*)$ can be naturally
viewed as either $\Cal G\x_P H_\ell(\frak p_+,\frak g)$ or $\Cal
G_0\x_{G_0} H_\ell(\frak p_+,\frak g)$. The latter interpretation
shows that they can be directly interpreted in terms of the underlying
structure. It is a crucial point in the theory that the
$G_0$--representations $H_\ell(\frak p_+,\frak g)$ can be computed
explicitly and algorithmically using Kostant's version of the
Bott--Borel--Weil theorem, see \cite{Yamaguchi}. (In that reference,
as well as in large parts of the literature, cohomology groups rather
than homology groups are used, but switching between the two points of
view is easy.)

A parabolic geometry $(p:\Cal G\to M,\om)$ is called \textit{normal}
if and only if its curvature $\ka$ has the property that
$\partial^*(\ka)=0$. 

\begin{thm*}
Let $(M,\{T^iM\})$ be a filtered manifold such that each symbol
algebra is isomorphic to $\fg_-$, and let $\Cal G_0\to M$ be a
reduction of $\gr(TM)$ to the structure group
$G_0\subset\Aut_{\gr}(\fg_-)$. Then there is a regular normal
parabolic geometry $(p:\Cal G\to M,\om)$ inducing the given data. If
$H_1(\frak p_+,\fg)$ is concentrated in non--positive homogeneous
degrees, then the pair $(\Cal G,\om)$ is unique up to isomorphism.
\end{thm*}

\subsection*{Remark} 
(1) The condition on $H_1(\frak p_+,\fg)$ can be easily turned into
something much more concrete, see \cite{Yamaguchi, Cap-Schichl}. If
$\frak g$ is simple, then it excludes exactly two series of examples
corresponding to the crossed Dynkin diagrams $\xbdbb$ and $\Cxbdbb$.
Except for the very degenerate case of the Dynkin diagram $\x$
(i.e.~the Borel subalgebra in $\frak{sl}(2,\Bbb K)$), the
corresponding regular normal parabolic geometries are still determined
by some underlying structure. Geometrically, these give rise to
classical projective structures and a contact version of projective
structures.

\noindent
(2) One actually obtains an equivalence of categories between regular
    normal parabolic geometries and regular infinitesimal flag
    structures.  

\subsection{Examples}\label{2.8}
By Theorem \ref{2.7}, a regular normal parabolic geometry on $M$ of
type $(G,P)$ is for almost all choices of $G$ and $P$ equivalent to a
filtration $\{T^iM\}$ of the tangent bundle such that each symbol
algebra is isomorphic to $\frak g_-$ plus a reduction of the structure
group of $\gr(TM)$ to the group $G_0$. In many situation, this
simplifies further, and we will discuss this next.

\medskip

\noindent
(1) \textbf{$|1|$--gradings.} 
Here we are in the situation $\frak g=\fg_{-1}\oplus\fg_0\oplus\fg_1$
and $\frak p=\frak g_0\oplus\frak g_1$. The classification of such
gradings is equivalent to the classification of Hermitian and
pseudo--Hermitian symmetric spaces and therefore well
known. Geometrically, the main point is that the filtration
$\{T^iM\}$ by definition consists of just one bundle. Moreover, the
regularity condition is easily seen to be vacuous in this case.

Hence if $(G,P)$ corresponds to a $|1|$--grading, then Theorem
\ref{2.7} says that normal parabolic geometries of type $(G,P)$ are
equivalent to classical first order $G_0$--structures on $M$. Here
$G_0$ is considered as a (covering of a) subgroup of $GL(\dim(M),\Bbb
R)$ via $\Ad:G_0\to GL(\fg_{-1})$. 

The most important examples of these structures are conformal, almost
quaternionic, and almost Grassmannian structures. The exceptional case
corresponding to the Dynkin diagram $\xbdbb$ corresponds to a
$|1|$--grading. Here $G_0=GL(\fg_{-1})$ so the underlying
infinitesimal flag structure contains no information at all. Normal
parabolic geometries of this type are equivalent to classical
projective structures, which will be discussed in more detail in
\ref{3.2} below. 

\medskip

\noindent
(2) \textbf{Structures determined by the filtration.}
We have seen in \ref{2.7} that the adjoint action defines a
homomorphism $G_0\to\Aut_{gr}(\fg_-)$. If this is an isomorphism, then
$\Cal G_0$ is the full natural frame bundle of $\gr(TM)$ and there is
no additional reduction of structure group. Hence in this case Theorem
\ref{2.7} shows that a regular normal parabolic geometry on $M$ is
equivalent to a filtration $\{T^iM\}$ such that each symbol algebra is
isomorphic to $\fg_-$. 

There is a simple way to obtain structures of this type: For any
semisimple $\fg$, the group $G:=\Aut(\fg)$ has Lie algebra $\fg$. It
turns out (see \cite{Katja}) that for this choice of $G$ we obtain
$G_0\cong\Aut_{gr}(\fg_-)$ provided that $H_1(\frak p_+,\fg)$ is
concentrated in negative homogeneous degrees. Again this homological
condition is easy to verify, and it turns out that it is often
satisfied. The paper \cite{Yamaguchi} contains a complete list of
pairs $(\fg,\fp)$ such that the condition is not satisfied.

This class of examples contains the quaternionic contact structures
introduced by O.~Biquard, see \cite{Biq, Biq2}, generic distributions of
rank 2 in dimension 5 (which were studied in Cartan's classic
\cite{Cartan:five}), rank 3 in dimension 6, and rank 4 in dimension 7.

\medskip

\noindent

(3) \textbf{Parabolic contact structures.}  These correspond to
$|2|$--gradings such that $\fg_{-2}$ is one--dimensional and such that
the bilinear form $\fg_{-1}\x\fg_{-1}\to\fg_{-2}$ defined by the
bracket is non degenerate. The classification of such gradings is
equivalent to the classification of quaternionic symmetric spaces and
therefore well know. Gradings of this type exist only on simple Lie
algebras and are unique up to isomorphism. With a few exceptions, they
exist on all non--compact, non--complex simple Lie algebras. 

Since $\fg_-$ by definition is a real Heisenberg algebra, a filtration
$TM=T^{-2}M\supset T^{-1}M$ of $TM$ such that each symbol algebra is
isomorphic to $\fg_-$ is exactly a contact structure $T^{-1}M\subset
TM$. Hence the filtration cannot be enough to determine the geometry
and one needs the additional reduction to the structure group $G_0$,
which can be expressed as an additional structure on $T^{-1}M$. 

This class contains non--degenerate partially integrable almost CR
structures of hypersurface type, for which the additional structure on
$T^{-1}M$ is an almost complex structure, as well as Lagrangean
contact structures, where the additional structure is a decomposition
of $T^{-1}M$ into the direct sum of two isotropic subbundles. Next,
there is the example of Lie contact structures (see
\cite{Sato-Yamaguchi}), in which the additional structure is a
decomposition of $T^{-1}M$ as the tensor product of two auxiliary
bundles, one of which has rank $2$ while the other one is endowed with
a pseudo--euclidean metric of some fixed signature. Finally, this
class also contains the second exceptional structure mentioned in
Remark \ref{2.7} (2). In that case, regular normal parabolic
geometries are equivalent to a contact analog of projective
structures, see \cite{Fox}.

\medskip

\noindent
(4) As an example of general parabolic geometries, we discuss
\textit{generalized path geometries}. These correspond to the
$|2|$--grading on $\frak{sl}(n+2,\Bbb R)$ corresponding to the first
and second simple root. In block form, this decomposition has the form
$$
\begin{pmatrix}
\fg_0 & \fg_1^L & \fg_2\\
\fg_{-1}^L & \fg_0 & \fg_1^R\\
\fg_{-2} & \fg_{-1}^R & \fg_0,
\end{pmatrix}
$$
where the blocks are of size $1$, $1$, and $n$. We have met this
grading for $n=2$ in \ref{2.4}. For later use, we have indicated
decomposition of $\fg_{\pm 1}$ into a one--dimensional part $\fg_{\pm
  1}^L$ and an $n$--dimensional part $\fg_{\pm 1}^R$.  Evidently, this
decomposition is invariant under the adjoint action of $\fg_0$. For an
appropriate choice of $G$, the subgroup $G_0$ consists of all
automorphisms of the graded Lie algebra $\fg_-$ which in addition
preserve the decomposition $\fg_{-1}=\fg_{-1}^L\oplus\fg_{-1}^R$.

From this description, we can directly read off the geometric meaning
of a regular infinitesimal flag structure of type $(G,P)$ on a smooth
manifold $M$ of dimension $2n+1$: One has two transversal subbundles
$L,R\subset TM$ of rank $1$ and $n$, respectively, such that for
$\xi,\eta\in\Ga(R)$ we have $[\xi,\eta]\in\Ga(L\oplus R)$ while the
Lie bracket induces an isomorphism $L\otimes R\to TM/(L\oplus R)$.

Examples of such structures come from path geometries. Let $N$ be a
manifold of dimension $n+1$ and consider the projectivized tangent
bundle $M:=\Cal PTN$, the space of lines through the origin in $TN$.
Take $R$ to be the vertical bundle of the projection $\Cal PTN\to N$.
Since $M$ is a projectivized tangent bundle, there is a tautological
subbundle $H\subset TM$ of rank $n+1$. The fiber of $H$ in a point
consists of those tangent vectors whose image in $TN$ lies in the line
determined by the point. Hence $R$ is contained in $H$ and a
\textit{path geometry} on $N$ is given by the choice of a line
subbundle $L\subset H$ such that $H=L\oplus R$. A path geometry on $N$
is equivalent to a family of unparametrized curves in $N$, with
exactly one curve through each point in each direction. In particular,
a system of second order ODE's on a manifold $Y$ can be equivalently
described as a path geometry on $Y\x\Bbb R$ by considering the
unparametrized curves describing the graphs of solutions, see
\cite{Grossman,Fels}.

For $n\neq 2$, the data $(M,L,R)$ corresponding to a regular
infinitesimal flag structure as above are locally isomorphic to a path
geometry. Namely, for $n\neq 2$ the subbundle $R\subset TM$ turns out
to be automatically integrable, and one defines $N$ to be a local leaf
space for the corresponding foliation. Then for an open subset
$U\subset M$, there is a surjective submersion $\ps:U\to N$ such that
$\ker(T_x\ps)=R_x$ for all $x\in U$.  Under $T_x\ps$, the line $L_x$
gives rise to a line in $T_{\ps(x)}N$, hence defining a lift
$\tilde\ps:U\to\Cal PTN$.  Possibly shrinking $U$, $\tilde\ps$ is an
open embedding. By construction, $T\tilde\ps$ maps $R$ to the vertical
subbundle and $L\oplus R$ to the tautological subbundle.

\subsection{Harmonic curvature}\label{2.9} 
There is a last element of the general theory of parabolic geometries
that we have to discuss. The Cartan curvature $\ka\in\Om^2(M,\Cal AM)$
as defined in \ref{2.2} and \ref{2.7} is a fairly complicated
object. In particular, to understand it geometrically, one needs the
adjoint tractor bundle, which is an equivalent encoding of the
principal Cartan bundle. An important feature of regular normal
parabolic geometries is that one may pass to the harmonic curvature
$\ka_H$, which is much easier to handle, but as powerful as $\ka$. 

In \ref{2.7} we have defined the operators $\partial^*:\La^\ell
T^*M\otimes\Cal AM\to \La^{\ell-1}T^*M\otimes\Cal AM$ and noted that
$\partial^*\o\partial^*=0$. For a normal geometry the curvature $\ka$
by definition is a section of the subbundle
$\ker(\partial^*)\subset\La^2T^*M\otimes\Cal AM$. Hence we can project
it to a section $\ka_H$ of the quotient
$\ker(\partial^*)/\im(\partial^*)$. As we have noted in \ref{2.7},
this quotient bundle can be identified with $\Cal G_0\x_{G_0}
H_2(\frak p_+,\fg)$, so it admits a direct interpretation in terms of
the underlying structure and is algorithmically computable.

We have also seen that $H_2(\frak p_+,\fg)$ splits into a direct sum
of $G_0$--irreducible components. Correspondingly, we obtain a
splitting of $\ka_H$ into fundamental curvature quantities. There are
several general tools to describe (parts of) $\ka_H$ in terms of the
underlying structure. 

The following result shows that $\ka_H$ still is a complete
obstruction to local flatness, and indeed, it contains the full
information about $\ka$. 
\begin{thm*}
Let $(p:\Cal G\to M,\om)$ be a regular normal parabolic geometry of
type $(G,P)$ with curvature $\ka$ and harmonic curvature $\ka_H$.

\noindent
(1) (Tanaka) If $\ka_H$ vanishes identically, then $\ka$ vanishes
    identically. 

\noindent
(2) (Calderbank--Diemer) There is a natural linear differential
operator $L$ such that $L(\ka_H)=\ka$.
\end{thm*}

The first part is a rather easy application of the Bianchi--identity
for Cartan connections. The second part is much more difficult. It
follows from the general machinery of BGG--sequences, see
\cite{CSS-BGG,David-Tammo}. 

\section{Correspondence spaces and twistor spaces}\label{3}

Now we switch to the discussion of constructions relating parabolic
geometries of different type. We start with the constructions of
correspondence spaces and twistor spaces, which is related to
different parabolic subgroups of the same group $G$. The basic
reference for this chapter is \cite{Cap:tw}. 

\subsection{Correspondence spaces}\label{3.1}
Consider a semisimple Lie group $G$ with nested parabolic subgroups
$Q\subset P\subset G$. For the homogeneous models, we have the simple
observation that $G/Q$ naturally fibers over $G/P$. Moreover, we can
interpret $G/Q$ as $G\x_P (P/Q)$, so this is the total space of a
natural fiber bundle over $G/P$. It turns out that the fiber $P/Q$ can
be equivalently viewed as the quotient of the semisimple part of
$G_0\subset P$ by its intersection with $Q$. This intersection turns
out to be parabolic, so $P/Q$ again is a generalized flag manifold.
The situations covered by this constructions are easy to describe in
the Dynkin (or Satake) diagram notation: The diagram corresponding to
$\frak q$ is obtained from the one corresponding to $\frak p$ by
changing dots into crosses. The fiber $P/Q$ can then be directly read
off the two diagrams, see \cite{B-E}.

Carrying this over to curved Cartan geometries is easy. Given a
geometry $(p:\Cal G\to N,\om)$ of type $(G,P)$ the subgroup $Q\subset
P$ acts freely on $\Cal G$. Hence the \textit{correspondence space}
$\Cal CN:=\Cal G/Q$ is a smooth manifold, and the obvious map $\Cal
G\to\Cal CN$ is a $Q$--principal bundle. Moreover, $\Cal CN=\Cal
G\x_P(P/Q)$, so $\pi:\Cal CN\to N$ is a natural fiber bundle with
fiber a generalized flag manifold. In particular, this fiber is always
compact. By definition, $\om\in\Om^1(\Cal G,\frak g)$ can also be
viewed as a Cartan connection on the principal $Q$--bundle $\Cal
G\to\Cal CN$.

The next obvious question is whether this construction is compatible
with regularity and normality. At this point, the uniform algebraic
construction of the normalization condition pays off:
\begin{prop*}
If $(\Cal G\to N,\om)$ is a normal parabolic geometry of type
$(G,P)$ then the parabolic geometry $(\Cal G\to\Cal CN,\om)$ of
type $(G,Q)$ is normal, too.
\end{prop*}
As we shall see in an example below, regularity is not preserved by
the construction in general. However, finding conditions which are
equivalent to regularity is usually very easy. 

\subsection{Example}\label{3.2}
Let $Q\subset G:=SL(n+2,\Bbb R)$ be the parabolic subgroup
corresponding to generalized path geometries as in Example (4) of
\ref{2.8}. Then $Q$ is the stabilizer of the flag consisting of the
line spanned by the first vector sitting inside the plane spanned by
the first two vectors of the standard basis of $\Bbb R^{n+2}$. Hence
we can write it as the intersection $P_1\cap P_2$ for parabolics $P_1$
and $P_2$ (the stabilizers of the line respectively the plane). Let us
start by analyzing the nested parabolics $Q\subset P_1\subset G$. 

Parabolic geometries of type $(G,P_1)$ correspond to classical
projective structures on $(n+1)$--dimensional manifolds, see Example
(1) of \ref{2.8}. Such a structure on a manifold $Z$ is given by the
choice of a projective equivalence class $[\nabla]$ of torsion free
linear connections on $TZ$. Two linear connections $\nabla$ and
$\hat\nabla$ on $TZ$ are called \textit{projectively equivalent} if
there is a one form $\Upsilon\in\Om^1(Z)$ such that
$$
\hat\nabla_\xi\eta=\nabla_\xi\eta+\Upsilon(\xi)\eta+\Upsilon(\eta)\xi
$$
for all vector fields $\xi,\eta\in\frak X(Z)$. Evidently,
projectively equivalent connections have the same torsion.
Alternatively, projective equivalence can be characterized as having
the same torsion and the same geodesics up to parametrization. The
harmonic curvature for this geometry is the projective Weyl curvature,
i.e.~the totally tracefree part of the curvature of any connection in
the class.

Since $\om$ is a Cartan connection on $\Cal G\to Z$, we have $TZ=\Cal
G\x_{P_1}(\frak g/\frak p_1)$. One easily verifies that $Q\subset P_1$
can be described as the stabilizer of a line in $\frak g/\frak p_1$.
Since $P_1$ acts transitively on the projective space $\Cal P(\frak
g/\frak p_1)$, see that $P/Q\cong\Cal P(\frak g/\frak p_1)$. Hence
$\Cal CZ=\Cal G\x_{P_1}P/Q$ can be naturally identified with the
projectivized tangent bundle $\Cal PTZ$. Since projective structures
are torsion free, the curvature $\ka$ of $\om$ has values in $\frak
p_1$, which immediately implies that $\om$ is regular as a Cartan
connection on $\Cal G\to\Cal CZ$. From Example (4) of \ref{2.8} we
conclude that $(\Cal G\to\Cal CZ,\om)$ can be interpreted as a path
geometry on $Z$. One verifies that the paths described in that way are
exactly the unparametrized geodesics of the connections from the
projective class.

Let us now switch to the nested parabolic subgroups $Q\subset
P_2\subset G$. A normal parabolic geometry $(\Cal G\to N,\om)$ of type
$(G,P_2)$ exists only for $\dim(N)=2n$ and is equivalent to an almost
Grassmannian structure. Essentially, such a geometry is given by two
auxiliary vector bundles $E$ and $F$ over $N$ of rank $2$ and $n$,
respectively, and an isomorphism $E\otimes F\to TN$. The subgroup
$Q\subset P_2$ can be characterized as the stabilizer of a line in the
representation inducing $E$, which similarly as above implies that
$\Cal CN$ can be identified with the projectivization $\Cal PE$ of
$E\to N$.
  
Here $\om$ is not regular as a Cartan connection on $\Cal G\to\Cal CN$
in general. Regularity turns out to be equivalent to the fact that the
structure on $N$ is Grassmannian rather than almost Grassmannian. This
can be characterized by vanishing of a certain torsion or equivalently
by the fact that there is a torsion free connection compatible with
the structure. If this is satisfied, then we obtain a generalized path
geometry on $\Cal PE$. The subbundle $L$ which is one of the
ingredients of that structure is simply the vertical bundle of $\Cal
PE\to N$. In particular, the manifold $N$ can be viewed the space of
all paths of the induced path geometry. The subbundle $L\oplus
R\subset T\Cal CN$ is again a tautological subbundle. The splitting of
this tautological subbundle as $L\oplus R$ comes from the torsion free
connections compatible with the Grassmannian structure. 

Suppose that $n>2$ (the case $n=2$ will be discussed later). Then
starting from a Grassmannian structure on $N$, we obtain a generalized
path geometry on $\Cal CN:=\Cal PE$. We know that the resulting
subbundle $R\subset T\Cal CN$ is involutive, so for sufficiently small
open subsets $U\subset\Cal CN$ we can form a local leaf space
$\ps:U\to Z$. With a bit more work, one shows that one may take
$U=\pi^{-1}(V)$, for sufficiently small and convex open subsets
$V\subset N$, where $\pi:\Cal CN\to N$ is the natural projection. One
then obtains a \textit{correspondence} 
$$
Z\overset{\ps}{\longleftarrow}
\pi^{-1}(V)\overset{\pi}{\longrightarrow} V,
$$
which is the basis for \textit{twistor theory} for Grassmannian
structures.

\subsection{Characterizing correspondence spaces}\label{3.3}
A central feature of the general theory of correspondence spaces is
that one can completely characterize parabolic geometries which are
locally isomorphic to correspondence spaces. This characterization is
uniform for all the structures.

Let us return to the general setting of nested parabolics $P\subset
Q\subset G$. The question we want to address is when a regular normal
parabolic geometry $(p:\Cal G\to M,\om)$ of type $(G,Q)$ is locally
isomorphic to the correspondence space $\Cal CN$ for a parabolic
geometry of type $(G,P)$. There is a fairly obvious necessary
condition: The subspace $\frak p/\frak q\subset\frak g/\frak q$ is
$Q$--invariant, thus giving rise to a subbundle $\Cal V\subset
TM$. For a correspondence space $\Cal CN$, this subbundle is the
vertical subbundle of the natural projection $\Cal CN\to N$. Since the
Cartan connections for $N$ and $\Cal CN$ are the same, so are their
curvatures. Since vectors from $\Cal V$ are vertical from the point of
view of $N$, they must hook trivially into the Cartan curvature of
$\Cal CN$.

It turns out that this condition is also sufficient:
\begin{thm*}
Let $(p:\Cal G\to M,\om)$ be a parabolic geometry of type $(G,Q)$ with
Cartan curvature $\ka$, and let $\Cal V\subset TM$ be the distribution
corresponding to $\frak p/\frak q\subset\frak g/\frak q$. Then $M$
admits an open covering $\{U_i\}$ such that the restriction of $(\Cal
G\to M,\om)$ to each $U_i$ is isomorphic to the correspondence space
of some parabolic geometry of type $(G,P)$ if and only if $i_\xi\ka=0$
for all $\xi\in\Cal V$.
\end{thm*}

The proof of this theorem is not specifically ``parabolic'' and uses
only principal bundle geometry. One first shows that the curvature
condition in the theorem implies that the distribution $\Cal V\subset
TM$ is involutive. Hence $\Cal V$ gives rise to a foliation of $M$,
and one considers a local leaf space for this foliation, i.e.~an open
subset $U\subset M$ together with a surjective submersion $\ps:U\to N$
such that $\ker(T_x\ps)=\Cal V_x$ for all $x\in U$. For sufficiently
small $U$, one next constructs a diffeomorphism from an open subset of
$p^{-1}(U)\subset\Cal G$ onto an open subset of the trivial principal
bundle $N\x P\to N$, which satisfies a certain equivariancy
condition. This diffeomorphism is then used to carry over $\om$ to
this open subset of $N\x P$, and one proves that the resulting form
uniquely extends to all of $N\x P$ by equivariancy. It is easy to see
that this not only gives a parabolic geometry of type $(G,P)$ on $N$
but also an isomorphism (of parabolic geometries) between $U$ and an
open subset of $\Cal CN$.  

While this result is very satisfactory from a conceptual point of
view, it is difficult to apply in concrete cases, since the Cartan
curvature is a complicated object. From part (2) of Theorem \ref{2.9} we
know that for regular normal geometries there is a natural
differential operator $L$ which computes the Cartan curvature from the
harmonic curvature $\ka_H$, which is much easier to handle. This
operator is constructed using the machinery of BGG sequences and the
construction is explicit enough to lead to relations between algebraic
properties of $\ka$ and $\ka_H$.
\begin{prop*}
Let $(\Cal G\to M,\om)$ be a regular normal parabolic geometry of type
$(G,Q)$ with Cartan curvature $\ka$ and harmonic curvature $\ka_H$,
and let $\Cal V\subset TM$ be as above. If $i_\xi\ka_H=0$ for all
$\xi\in\Cal V$, then $i_\xi\ka=0$ for all $\xi\in\Cal V$.
\end{prop*}
Combining this result with the theorem above, one obtains a very
efficient local characterization of correspondence
spaces. From another point of view, these are equivalent conditions
for the existence of natural geometric structures on twistor
spaces. It has to be pointed out here that usually the structure of
the harmonic curvature can be understood without detailed knowledge of
the canonical Cartan connection. 

\subsection{Examples}\label{3.4}
Let us interpret the results on local characterization of
correspondence spaces in the example discussed in \ref{3.2}. So we
start with a generalized path geometry $(M,L,R)$ and the associated
regular normal parabolic geometry $(p:\Cal G\to M,\om)$ of type
$(G,Q)$. For $n>2$ (which we will still assume throughout this
subsection), the harmonic curvature $\ka_H$ splits into two
irreducible components:
\begin{align*}
    &T:L\wedge TM/(L\oplus R)\to R \qquad\text{Torsion}\\
    &\rho:R\wedge TM/(L\oplus R)\to R^*\otimes R\qquad\text{Curvature}
\end{align*}
The types of these components can be deduced from the structure of the
homology group $H_2(\frak q_+,\frak g)$, which can be determined
algorithmically using Kostant's version of the Bott--Borel--Weil
theorem. There are general procedures how to obtain explicit formulae
for the two components, say in terms of a local non--vanishing section
of $L$. 

Let us first consider the characterization of correspondence spaces
coming from the inclusion $Q\subset P_1\subset G$. From \ref{3.2} we
know that these are exactly the path geometries associated to the
unparametrized geodesics of a projective class of connections. The
distribution $\Cal V$ corresponding to $\frak p_1/\frak q\subset\frak
g/\frak q$ evidently is the subbundle $R\subset TM$. The results from
\ref{3.3} now show that $M$ is locally isomorphic to a correspondence
spaces if and only if $\rho$ vanishes identically. 

As we have noted in \ref{2.8}, the subbundle $R\subset TM$ is
involutive (since $n>2$). For a local leaf space $\ps:U\to Z$ of the
corresponding foliation, the subset $U$ then is naturally
diffeomorphic to an open subset in the projectivized tangent bundle
$\Cal PTZ$. Then our result shows that the generalized path geometry
on $M$ induces a projective structure on $Z$ if and only if $\rho$
vanishes identically. If this is the case, then the torsion $T$ is
directly related to the projective Weyl curvature of the induced
structures on the local leaf spaces. In particular, the path geometry
on $M$ is locally flat if and only if the induced projective
structures on all local leaf spaces are locally projectively flat.

Another interesting application of this criterion is to the path
geometry associated to a system of second order ODE's as described in
\ref{3.2}. This reproduces a result of \cite{Fels}:
\begin{thm*}
A system of second order ODE's is locally equivalent to a geodesic
equation if and only if the curvature $\rho$ of the associated path
geometry vanishes identically.
\end{thm*}

Now we switch to the characterization of correspondence spaces with
respect to the inclusion $Q\subset P_2\subset G$. The distribution
$\Cal V$ corresponding to $\frak p_2/\frak q\subset\frak g/\frak q$ is
the subbundle $L\subset TM$. This is always involutive and local leaf
spaces for the associated foliation locally parametrize the paths of
the path geometry. Hence here the main interpretation of the
characterization result is a criterion when a generalized path
geometry locally descends to a Grassmannian structure on the space of
all paths. From \ref{3.3} we see that this is the case if and only if
$T=0$, which is equivalent to the generalized path geometry being
torsion free. 

Again there is an interesting application to the theory of systems of
second order ODE's: One defines such a system to be torsion free if
and only if the associated path geometry is torsion free. For such a
systems we obtain an induced Grassmannian structure on the space of
solutions of the system. The curvature of this Grassmannian structure
can be constructed from the curvature $\rho$ of the path geometry. Of
course, this curvature descends to the space of solutions and hence is
constant along each solution. Using this, D.~Grossman proved in
\cite{Grossman} the following result.
\begin{thm*}
For generic torsion free systems of second order ODE's, the curvature
of the associated path geometry can be used to solve the system
explicitly.
\end{thm*}

\subsection{The case $n=2$}\label{3.5}
Let us briefly discuss how the examples related to generalized path
geometries discussed in \ref{2.9}, \ref{3.2}, and \ref{3.4} change in
the case $n=2$. The ingredients are projective structures on three
manifolds, generalized path geometries in dimension five, and four
dimensional almost Grassmannian structures. The main point is that an
almost Grassmannian structure in dimension four is equivalent to a
conformal pseudo--Riemannian spin structure of split signature
$(2,2)$. The auxiliary bundles $E$ and $F$ whose tensor product is
isomorphic to the tangent bundle both have rank two. They are exactly
the two spinor bundles. 

The structure of harmonic curvatures for $n=2$ is also different from
the case $n>2$. For almost Grassmannian structures the more symmetric
situation leads to the fact that there are two curvatures rather than
one curvature and one torsion. These two components are exactly the
self dual and the anti self dual part of the Weyl curvature of the
corresponding conformal structure. 

On the level of path geometries, a third irreducible component in the
harmonic curvature shows up. This component is represented by a
torsion $\tau:\La^2R\to L$, which is the obstruction to involutivity
of the subbundle $R$. (For $n>2$, there also is a corresponding
component in the homology $H_2(\frak q_+,\frak g)$, but this sits in
homogeneity zero. By regularity, this component cannot contribute to
the harmonic curvature.)

Starting from a conformal four manifold, the correspondence space is a
projectivized spinor bundle, which inherits a generalized path
geometry. The torsion $\tau$ on this space corresponds exactly to  the
self dual part of the Weyl curvature downstairs. Vanishing of this
part, i.e.~anti self duality, is equivalent to existence of local
leaf spaces for the bundle $R$ on the correspondence space. This is
the basis for twistor theory for anti self dual four manifolds in
split signature. The Riemannian version of twistor theory can be
either obtained from the complex version of this construction or by an
analog of the correspondence space construction (for a subgroup which
is not parabolic). 

\section{Analogs of the Fefferman construction}\label{4}
We now switch to a second general construction relating parabolic
geometries of different types. The basic example for this is
Fefferman's construction which relates CR structures to conformal
structures. This construction is of different nature to the ones
discussed in section \ref{3} since it involves two different
semisimple groups. More details on the contents of this section can be
found in \cite{Cap:Fefferman} and \cite{Cap-Gover:Fefferman}. 

\subsection{The Fefferman construction}\label{4.1}
We start by reviewing Fefferman's original construction from 
\cite{Fefferman} and its interpretation in terms of Cartan
geometries. He started from a strictly pseudoconvex domain
$\Om\subset\Bbb C^{n+1}$ with smooth boundary $M:=\partial\Om$. This
boundary naturally inherits a CR structure (see below). Studying the
Bergman kernel of $\Om$, Fefferman was led to consider the
\textit{ambient metric}: Put $\Bbb C^*:=\Bbb C\setminus\{0\}$ and
consider $M_{\#}=M\x\Bbb C^*\subset\Om_\#=\Om\x\Bbb C^*$. A defining
function $r$ for $M$ induces a defining function $r_\#$ for
$M_\#$. Since $M$ is strictly pseudoconvex, $r_\#$ can be used as the
potential for a pseudo--K\"ahler metric $g_\#$ of signature
$(n+1,1)$. Fefferman showed that one may always chose $r$ to be an
approximate solution of a Monge--Amp\`ere equation and doing this a
certain jet of $g_\#$ along $M_\#$ is invariant under biholomorphisms
of $\Om$. Otherwise put, this jet is a CR invariant of $M$.

Hence it is a natural idea to look at the restriction of $g_\#$ to
$M_\#$. This turns out to be degenerate but only in the real
directions within the vertical subspaces of the projection $M_\#\to
M$. To get rid of these directions, one passes to the space $\tilde
M=M\x(\Bbb C^*/\Bbb R^*)\cong M\x S^1$. Using a section of the evident
projection $M_\#\to\tilde M$, one can pull back $g_\#$ to a
non--degenerate Lorentz metric on $\tilde M$. Changing the sections
leads to a conformal change of the metric, so one obtains a well
defined conformal class of metrics of signature $(2n+1,1)$ on $\tilde
M$. This conformal class is invariant under biholomorphisms of $\Om$
and hence depends only on the CR structure of $M$.

CR structures fit into the general concept of parabolic geometries as
the parabolic contact structures associated to $\frak
g=\frak{su}(p+1,q+1)$. In fact, one obtains a more general concept: A
\textit{partially integrable almost CR structure} on a smooth manifold
$M$ of dimension $2n+1$ is a contact structure $H\subset TM$ together
with an almost complex structure $J:H\to H$ such that the Levi bracket
$\Cal L$ (see \ref{2.6}) satisfies $\Cal L(J\xi,J\eta)=\Cal
L(\xi,\eta)$ for all $\xi,\eta$. Under this assumption, $\Cal L$ is
the imaginary part of a Hermitian form (with values in the real line
bundle $TM/H$), the \textit{Levi form}, which has a signature $(p,q)$.
Since there is an ambiguity of sign, we require $p\geq q$ to have the
signature well defined.

The compatibility of $\Cal L$ and $J$, which is usually referred to as
partial integrability, can also be nicely formulated in terms of
complexifications. The almost complex structure $J$ leads to a
splitting of $H\otimes\Bbb C\subset TM\otimes\Bbb C$ into the direct sum
of the holomorphic part $H^{1,0}$ and the anti holomorphic part
$H^{0,1}$, which are conjugate to each other. Partial integrability is
equivalent to the fact that the Lie bracket of two sections of
$H^{0,1}$ is a section of $H\otimes\Bbb C$. An almost CR structure is
called \textit{integrable} or a \textit{CR structure} if the subbundle
$H^{0,1}M\subset TM\otimes\Bbb C$ is involutive, so the Lie bracket of
two sections of $H^{0,1}$ even is a section of $H^{0,1}$.

Partially integrable almost CR structures of signature $(p,q)$ are
then equivalent to regular normal parabolic geometries associated to
the group $PSU(p+1,q+1)$. For many applications it is better to extend
the group to $G:=SU(p+1,q+1)$. Let $P$ be the stabilizer of an
isotropic complex line $\ell$ in $\Bbb V:=\Bbb C^{p+q+2}$. Then a
regular normal parabolic geometry of type $(G,P)$ on a manifold $M$ is
equivalent to a partially integrable almost CR structure of signature
$(p,q)$ plus a choice of a complex line bundle, which is an $(n+2)$nd
root of the so--called canonical bundle. While such a choice need not
exist in general, it is always possible locally. The integrability
condition turns out to be equivalent to torsion freeness of the
associated parabolic geometry.

If $M$ is the boundary of a strictly pseudoconvex domain
$\Om\subset\Bbb C^{n+1}$, then one defines $H_xM$ as the maximal
complex subspace of $T_xM\subset T_x\Bbb C^{n+1}$. This evidently has
an almost complex structure and it defines a contact structure by
strict pseudoconvexity. The latter condition also implies that the
signature is $(n,0)$. Looking at the complexified tangent bundle, we
see that $H^{0,1}M=(TM\otimes\Bbb C)\cap T^{0,1}\Bbb C^{n+1}$. Since
$\Bbb C^{n+1}$ is a complex manifold, the subbundle $T^{0,1}\Bbb
C^{n+1}\subset T\Bbb C^{n+1}\otimes\Bbb C$ is involutive, so we obtain
a CR structure on $M$. Triviality of the tangent bundle of $\Bbb
C^{n+1}$ implies that the canonical bundle of $M$ is canonically
trivial, so there is no problem in choosing an $(n+2)$nd root. 

Now it is easy to obtain the Fefferman construction for the
homogeneous model: The real part of the Hermitian form on $\Bbb V$
defines an inner product of signature $(2p+2,2q+2)$ on the underlying
real vector space $\Bbb V_{\Bbb R}$. Since elements of $G$ preserve
this real part, we obtain an injection $G\hookrightarrow
SO(2p+2,2q+2)$. Analyzing the induced homomorphism between the
fundamental groups one even shows that this naturally lifts to an
inclusion into the spin group $\tilde G:=Spin(2p+2,2q+2)$. Choose a
real line $\ell_{\Bbb R}$ in the isotropic complex line $\ell$ and let
$\tilde P\subset\tilde G$ be the stabilizer of $\ell_{\Bbb R}$. The
intersection $Q:=G\cap\tilde P$ is the stabilizer of $\ell_{\Bbb R}$
in $G$, so it is evidently contained in $P$ and $P/Q\cong\Bbb
RP^1$. Elementary linear algebra shows that $G$ acts transitively on
the space of real null lines in $\Bbb V_{\Bbb R}$. Hence the inclusion
$G\hookrightarrow\tilde G$ induces a diffeomorphism $G/Q\to\tilde
G/\tilde P$. The latter space is well known to be the homogeneous
model of conformal spin structures of signature $(2p+1,2q+1)$. Hence
we obtain such a structure (which by construction is invariant under
the action of $G$) on $G/Q$ which is the total space of a circle
bundle over $G/P$.

Passing to curved geometries is easy: Looking at the tangent spaces at
the base points, the diffeomorphism $ G/Q\to\tilde G/\tilde P$ induces
a linear isomorphism $\frak g/\frak q\to \tg/\tp$ which is equivariant
over the inclusion $Q\hookrightarrow \tilde P$. Here $\frak q=\frak
g\cap\tp$ is the Lie algebra of $Q$. In particular, we obtain a
conformal class of inner products of signature $(2p+1,2q+1)$ on
$\frak g/(\frak g\cap\tp)$ which is invariant under the natural action
of $Q$. Given a partially integrable almost CR structure $(M,H,J)$,
let $(\Cal G\to M,\om)$ be the associated regular normal parabolic
geometry. The subgroup $Q\subset P$ acts freely on $\Cal G$, so the
\textit{Fefferman space} $\tilde M:=\Cal G/Q$ is a smooth manifold and
the total space of the natural fiber bundle $\Cal G\x_PP/Q$ over
$M$. On the other hand, the evident projection $\Cal G\to\tilde M$ is
a principal $Q$--bundle and $\om\in\Om^1(\Cal G,\frak p)$ defines a
Cartan connection on that bundle. In particular, $T\tilde M\cong \Cal
G\x_Q\frak g/\frak q$ so the $Q$--invariant class of inner products on
$\frak g/\frak q$ gives rise to a conformal structure on $\tilde M$,
which by construction depends only on the CR structure on $M$.

It is easy to give a more explicit description of $\tilde M$. Namely,
one shows that $\tilde M$ can be naturally identified with the space
of real lines in a natural complex line bundle, which is closely
related to the chosen root of the canonical bundle. One can also
construct explicitly a metric from the conformal class in terms of a
choice of contact form on $M$ (usually called a \textit{pseudo
Hermitian structure}) and the associated Weyl connection (see
\cite{Weyl}) on a complex line bundle. 

\subsection{Cartan geometry interpretation}\label{4.2}
The construction of the canonical conformal class on $\tilde M$ from
above can be easily reformulated in terms of Cartan geometries. As we
know from \ref{4.1}, we have the $Q$--principal bundle $\Cal G\to\tilde
M$ and we can view the canonical CR Cartan connection $\om$ as a
Cartan connection on that bundle. Now via the inclusion
$Q\hookrightarrow\tilde P$, we can extend the structure group of this
bundle. Define a principal $\tilde P$--bundle $\tcg:=\Cal G\x_Q\tilde
P\to\tilde M$. Mapping $u\in\Cal G$ to the class of $(u,e)$ in $\tcg$
defines an injective smooth map $j:\Cal G\to\tcg$ which is equivariant
over the inclusion $Q\hookrightarrow\tilde P$. It is easy to show that
there is a unique Cartan connection $\tilde\om\in\Om^1(\tcg,\tg)$ such
that $\tilde\om|_{Tj(T\Cal G)}=\om$ (viewing $\frak g$ as a Lie
subalgebra of $\tg$). 

As a Cartan connection on a principal $\tilde P$--bundle $\tilde\om$
is automatically regular and hence it induces a conformal spin
structure on the base $\tilde M$. From the construction it is evident
this this leads to the conformal structure described in \ref{4.1}. 

Now one might expect that $\tilde\om$ is the normal Cartan connection
associated to this conformal spin structure, but this is \textit{not}
true in general:
\begin{thm*}
Let $(M,H,J)$ be a partially integrable almost CR structure with 
Fefferman space $\tilde M$. Then the Cartan connection $\tilde\om$ on
the extended principal bundle $\tcg\to\tilde M$ is normal if and only
if the almost CR structure is integrable. 
\end{thm*}

The necessity of integrability follows rather easily from the fact
that normal conformal Cartan connections are automatically torsion
free. The proof of sufficiency of this condition is much more
subtle. The result does not follow from algebraically comparing the
normalization conditions for the two geometries in question but one
has to prove additional properties of the curvature of a torsion free
geometry. In that respect, the situation is very different from the
case of correspondence spaces discussed in the last section. 

For some applications of the Fefferman construction, the question of
normality of $\tilde\om$ is not relevant. For example, conformal
invariants of the Fefferman space are always invariants of the
underlying partially integrable almost CR structure. However, we will
show below that normality of $\tilde\om$ leads to many other and
deeper results. 

If the structure on $M$ is not integrable, then the canonical Cartan
connection for the conformal spin structure on $\tilde M$ can be
obtained by normalizing $\tilde\om$. The difference of $\tilde\om$
from the normal Cartan connection is given by a one form on $\tilde M$
with values in the conformal adjoint tractor bundle $\tcg\x_{\tilde
P}\tg$. One may try to imitate some of the developments described
below taking into account the change caused by this form. To my
knowledge, this has not been explored up to now.

\subsection{Applications of normality to CR geometry}\label{4.3} 
We want to discuss a few results which are based on normality of the
Cartan connection $\tilde\om$ in the case of a CR structure. The first
of these was the main application in Fefferman's original article
\cite{Fefferman} as well as in the first version for abstract CR
structures in \cite{BDS}.

\medskip

\noindent
$\bullet$ \textit{Chern--Moser chains are the projections to $M$ of
null geodesics in $\tilde M$}.

Chern--Moser chains in $M$ can be obtained as the projections of flow
lines of vector fields on $\Cal G$ which are mapped to certain
constant functions by $\om$. Likewise, conformal circles on $\tilde M$
are the projections of flow lines of vector fields on $\tcg$ which are
mapped to certain constant functions by $\tilde\om$. For initial
directions which are null, conformal circles are just null geodesics
which, as unparametrized curves, are conformally invariant. The initial
direction of a chain is always transversal to the contact subbundle,
and such a direction always admits a lift to a null direction in
$\tilde M$. Then the result easily follows from the fact that
$\tilde\om$ is obtained from $\om$ by equivariant extension. 

\medskip

\noindent
$\bullet$ \textit{Relations between CR tractor calculus on $M$ and
  conformal tractor calculus on $\tilde M$}.

Standard tractors are probably the nicest way to relate a CR manifold
to its Fefferman space. The CR standard tractor bundle $\Cal T$ of $M$
is by definition the associated bundle $\Cal G\x_P\Bbb V$, where $\Bbb
V$ denotes the standard representation of $G$. By construction, this
is a rank $n+2$ complex vector bundle endowed with Hermitian inner
product $h$ of signature $(p+1,q+1)$, and a complex line subbundle
$\Cal T^1\subset\Cal T$ which is isotropic for $h$. This subbundle
corresponds to the complex line in $\Bbb V$ which is stabilized by
$P$. The canonical Cartan connection $\om$ on $\Cal G$ induces a
Hermitian linear connection on $\Cal T$, called the \textit{normal
standard tractor connection}. 

Likewise, the conformal standard tractor bundle $\tct$ of the
Fefferman space $\tilde M$ is the bundle $\tcg\x_{\tilde P}\Bbb
V$. This is a real bundle of rank $2n+4$ endowed with a Euclidean
bundle metric $\tilde h$ of signature $(2p+2,2q+2)$ and an a real line
subbundle $\tct^1$ which is isotropic for $\tilde h$. The Cartan
connection $\tilde\om$ induces the normal standard tractor connection
on $\tct$. 

The relation between the Cartan bundles and the Cartan connections
discussed above can be interpreted as the fact that $\tct$ (including
the additional structures) can also be obtained as $\Cal
G\x_{G\cap\tilde P}\Bbb V$ and the normal tractor connection on $\tct$
is induced by $\om$, viewed as a Cartan connection on $\Cal G\to\tilde
M$. 

Both for conformal and for CR structures, the standard tractor bundle
and the standard tractor connection lead to an efficient
calculus. Hence we obtain a close relation between CR tractor calculus
on a CR manifold and conformal tractor calculus on its Fefferman
space. 

\medskip

\noindent
$\bullet$ \textit{Conformally invariant differential operators on
$\tilde M$ descend to families of CR invariant differential operators
on $M$.}

The relations between the standard tractor bundles of $M$ and $\tilde
M$ can be extended to other bundles, for example other tractor bundles
and density bundles. One can then interpret sections of some bundle
over $M$ as a subset of sections of some other bundle over $\tilde M$,
which usually are characterized as solutions of some differential
equation. It often happens that this works for a whole family of
bundles over $M$ (with different weights) and the same bundle on
$\tilde M$. Based on the relations between tractor calculi discussed
above, one shows that in several cases conformally invariant
differential operators preserve the subspaces of ``downstairs''
sections and hence descend to (families of) CR invariant differential
operators.

\medskip

\noindent
$\bullet$ \textit{Interpretation of solutions of certain CR
invariant differential equations.}

The solutions of certain CR invariant differential equations admit a
natural interpretation in terms of the conformal geometry of the
Fefferman space. An example for this will be given in the discussion
of conformal isometries of the Fefferman space below. 

\subsection{Conformal geometry of Fefferman spaces}\label{4.4}
The second interesting line of applications is towards Fefferman
spaces as an interesting subclass of conformal structures. 

\medskip

\noindent
$\bullet$ \textit{Fefferman spaces have a parallel orthogonal complex
structure on the standard tractor bundle and are locally characterized
by that.}

We have seen above that for the Fefferman space $\tilde M$ of a CR
manifold $M$, the conformal standard tractor bundle can be interpreted
as $\tct=\Cal G\x_{G\cap\tilde P}\Bbb V$, and the tractor connection
on that bundle is induced by the CR Cartan connection $\om$. Since
$\Bbb V$ is a complex vector space, we obtain an almost complex
structure $J$ on $\tct$, which is orthogonal (or equivalently skew
symmetric) with respect to the tractor metric, and parallel for the
connection on $L(\tct,\tct)$ induced by the standard tractor
connection. 

This can be interpreted as the fact that the holonomy of the standard
tractor connection is contained in $SU(p+1,q+1)\subset
SO(2p+2,2q+2)$. Conversely, one can show that a conformal structure of
signature $(2p+2,2q+2)$ which admits such a holonomy reduction, is
locally conformally isometric to a Fefferman space. This shows that
the role of Fefferman spaces among general conformal structures is
similar to the role of Calabi--Yau manifolds among general Riemannian
manifolds. 

\medskip

\noindent
$\bullet$ \textit{Fefferman spaces admit nontrivial Twistor spinors
  and conformal Killing forms of all odd degrees.}

Several conformally invariant differential equations which are
overdetermined (and thus do not have solutions in general) always
admit nontrivial solutions on Fefferman spaces. The simplest example of
this situation is that one constructs a nowhere vanishing conformal
Killing field $j$ on $\tilde M$, which spans the vertical subbundle of
$\tilde M\to M$. The most conceptual interpretation of this is via the
almost complex structure $J$ on the standard tractor bundle
$\tct\to\tilde M$. Since this is skew symmetric with respect to the
tractor metric, it can be interpreted as a parallel section of the
adjoint tractor bundle $\tca=\tcg\x_{\tilde P}\tg$. It is well known
that there is a natural projection $\Pi:\tca\to T\tilde M$ and the
image of a parallel section under this projection is automatically a
conformal Killing field (which in addition hooks trivially into the
Cartan curvature). 

Viewing $\tca$ as $\La^2\tct$, we can form the $k$--fold wedge product
of $J$ with itself, which defines a nonzero parallel section of the
tractor bundle $\La^{2k}\tct$. This bundle naturally projects onto the
bundle $\La^{2k-1}T^*\tilde M$ (twisted by an appropriate density
bundle) and the image of a parallel section is a conformal Killing
form (with additional properties), see \cite{Leitner}. These conformal
Killing forms can be explicitly expressed in terms of the conformal
Killing field $j$ from above. In contrast to the simple algebraic
formula on the tractor level, these expressions involve covariant
derivatives of $j$.

We have noted in \ref{4.1}, the Fefferman space $\tilde M$ carries a
natural spin structure. In particular, we can consider the tractor
bundle $\tcs\to\tilde M$ corresponding to the spin representation of
$\tilde G=Spin(2p+2,2q+2)$. Now it is well known that as a
representation of the subgroup $G=SU(p+1,q+1)$ this spin
representation decomposes and in particular contains a two dimensional
trivial subrepresentation. Using the relation between the tractor
calculi discussed above, one shows that this leads to a decomposition
of the spin tractor bundle $\tcs\to\tilde M$, and in particular one
obtains a two parameter family of parallel sections of that
bundle. The bundle $\tcs$ comes with a canonical projection to the
spinor bundle of $\tilde M$, which maps parallel sections to twistor
spinors. Hence any Fefferman space admits a two parameter family of
twistor spinors. 

\medskip

\noindent
$\bullet$ \textit{Decomposition of conformal Killing fields.}

By naturality of the construction of the Fefferman space, any
CR automorphism of $M$ lifts to a conformal isometry of $\tilde
M$. Likewise, an infinitesimal automorphism of $M$ induces a conformal
Killing field on $\tilde M$. As the example of the homogeneous model
shows, there may be other conformal Killing fields on $\tilde M$. It
turns out that one can completely describe the space of all conformal
Killing fields on $\tilde M$ in terms of the CR geometry of $M$. 

Infinitesimal automorphisms of parabolic geometries can be described in
general in terms of sections of the adjoint tractor bundle. For the
case of conformal structures, this means that any conformal Killing
field is the image of a uniquely determined section of the adjoint
tractor bundle $\tca$ which satisfies a certain conformally invariant
differential equation.

As a representation of $G=SU(\Bbb V)$, the Lie algebra
$\tg=\frak{so}(\Bbb V)$ is not irreducible, but decomposes as
$\frak{su}(\Bbb V)\oplus\Bbb R\oplus\La^2_{\Bbb C}\Bbb V$. Here the
first two summands correspond to complex linear maps, while the last
one corresponds to conjugate linear maps, and the trivial summand
consists of purely imaginary multiples of the identity. This induces an
analogous splitting of the conformal adjoint tractor bundle
$\tca\to\tilde M$.

We can use this splitting to decompose any section of $\tca$ into a
sum of three terms. Via tractor calculus one shows that for a section
corresponding to a conformal Killing field, each of the three parts
satisfies the infinitesimal automorphism equation. Thus one concludes
that any conformal Killing fields $\xi\in\frak X(\tilde M)$ decomposes
uniquely into a sum $\xi_1+\xi_2+\xi_3$ of conformal Killing fields.
One further shows that $\xi_1$ descends to an infinitesimal
automorphism of the underlying CR manifold $M$ and $\xi_2$ is a
constant multiple of $j$. The summand $\xi_3$ descends to a section of
$\La^2_{\Bbb C}\Cal T\to M$ which solves a certain CR invariant
differential equation. Likewise, appropriate solutions of this
equation give rise to conformal Killing fields on $\tilde M$.

\subsection{Analogs of the Fefferman construction}\label{4.5}
From the discussion in \ref{4.1} it is pretty evident what is needed
to obtain an analog of the Fefferman construction: One starts with an
inclusion $G\hookrightarrow\tilde G$ of semisimple Lie groups and
chooses a parabolic subgroup $\tilde P\subset\tilde G$ such that the
$G$ orbit of $e\tilde P$ in $\tilde G/\tilde P$ is open. Finally, one
needs a parabolic subgroup $P\subset G$ which contains $G\cap\tilde
P$.

Suppose that $(p:\Cal G\to M,\om)$ is a parabolic geometry of type
$(G,P)$. The define $\tilde M:=\Cal G/(G\cap\tilde P)$, which is a
smooth manifold and the total space of the natural fiber bundle $\Cal
G\x_PP/(G\cap\tilde P)\to M$. To obtain an explicit description of
$\tilde M$, it suffices to give a good description of the subgroup
$G\cap\tilde P\subset P$. As before, one can view $\Cal G\to\tilde M$
as a principal bundle with structure group $G\cap\tilde P$ and
$\om\in\Om^1(\Cal G,\frak g)$ as a Cartan connection on this
bundle. In particular, this identifies $T\tilde M$ with the associated
bundle $\Cal G\x_{G\cap\tilde P}\frak g/(\frak g\cap\tp)$. 

Since the $G$--orbit of $e\tilde P$ in $\tilde G/\tilde P$ is open,
the inclusion $\frak g\hookrightarrow\tg$ induces a linear isomorphism
$\frak g/(\frak g\cap\tp)\to\tg/\tp$. Clearly, this isomorphism is
equivariant under the inclusion $G\cap\tilde P\hookrightarrow\tilde
P$. Hence we can carry over $\tilde P$--invariant objects related to
$\tg/\tp$ to $(G\cap\tilde P)$--invariant objects related to $\frak
g/(\frak g\cap\tp)$ and hence to natural geometric objects on $\tilde
M$. In most examples discussed below, this already suffices to obtain
the underlying structure of a regular normal parabolic geometry of
type $(\tilde G,\tilde P)$ on $\tilde M$. In more complicated
situations one in addition has to check that the map
$\La^2(\tg/\tp)\to\tg$ induced by the curvature of $\om$ is regular,
but this usually is very easy. 

It is a much more difficult problem to check whether $\om$ induces the
regular normal Cartan connection associated to this underlying
structure.  As in the classical case, one can form the extended bundle
$\tcg:=\Cal G\x_{G\cap\tilde P}\tilde P$, and there is a unique Cartan
connection $\tilde\om$ on that bundle which restricts to $\om$ on
$T\Cal G\subset T\tcg$. To obtain an analog of Theorem \ref{4.2} and
applications similar to the ones described in \ref{4.3} and \ref{4.4},
one has to find conditions for $\tilde\om$ being normal.  To my
knowledge, this has not been done for all the examples described below
but for many of them there are hints coming from independent works on
these structures.

\subsection*{Examples}
(1) Closest to the classical Fefferman construction, one may consider
the group $G:=Sp(p+1,q+1)$ associated to a quaternionic Hermitian form
of signature $(p+1,q+1)$ on $\Bbb H^{p+q+2}$. Viewing this space as
$\Bbb C^{2p+2q+4}$ gives rise to an inclusion
$Sp(p+1,q+1)\hookrightarrow \tilde G:=SU(2p+2,2q+2)$. Taking $P\subset
G$ and $\tilde P\subset\tilde G$ the stabilizer of a quaternionic
respectively a complex null line, one obtains $G\cap\tilde P\subset P$
and $P/(G\cap\tilde P)\cong \Bbb CP^1$. Parabolic geometries of type
$(G,P)$ fall into the class discussed in Example (2) of \ref{2.8},
i.e.~the structures which are (essentially) determined by a filtration
of the tangent bundle. The modeling Lie algebra $\frak g_-$ is a
quaternionic Heisenberg algebra of signature $(p,q)$. This means that
$\frak g_{-1}\cong\Bbb H^{p+q}$ and $\frak g_{-2}\cong\Im(\Bbb H)$,
the space of purely imaginary quaternions, in such that way that the
bracket is by the imaginary part of a quaternionic Hermitian form of
signature $(p,q)$. For $q=0$, one obtains the quaternionic contact
structures introduced by Olivier Biquard, see \cite{Biq, Biq2}.

Hence we see that, up to some discrete data (related to the fact that
we use the group $Sp$ rather than $PSp$) our construction starts with
a quaternionic contact structure of signature $(p,q)$ on some manifold
$M$. The Fefferman space $\tilde M$ is then the total space of a
natural fiber bundle over $M$ with fiber $\Bbb CP^1\cong S^2$, and on
$\tilde M$ we naturally obtain a partially integrable almost CR
structure of signature $(2p+1,2q+1)$. This should be closely related
to O.~Biquard's construction of a twistor space for quaternionic
contact structures. 

\noindent
(2) Consider a vector space $\tilde{\Bbb V}$ endowed with an inner
product of signature $(p+1,q+2)$. Fixing a line $\ell$ on which the
inner product is negative definite, the inclusion
$\ell^\perp\hookrightarrow\tilde{\Bbb V}$ gives rise to an inclusion
$G:=SO(p+1,q+1)\hookrightarrow SO(p+1,q+2)=:\tilde G$. Choose a null
plane $\Bbb W$ which is transversal to $\ell^\perp$ and let $\tilde
P\subset\tilde G$ be the stabilizer of $\Bbb W$. Then $\Bbb
W\cap\ell^\perp$ is a null line, and its stabilizer $P$ evidently
contains $G\cap\tilde P$. One verifies that the $G$--orbit of $e\tilde
P$ in $\tilde G/\tilde P$ consists of those null planes in
$\tilde{\Bbb V}$ which are transversal to $\ell^\perp$, so in this
case $G/(G\cap\tilde P)$ is a proper open subset of $\tilde G/\tilde
P$. 

Normal parabolic geometries of type $(G,P)$ are just conformal
structures of signature $(p,q)$. Given such a structure on $M$ one
shows that the Fefferman space $\tilde M$ can be identified with the
open subset $\Cal P_+(T^*M)$ of the projectivized cotangent bundle of
$M$ consisting of all lines in $T^*M$ on which the conformal inner
product is positive definite. In particular, for $q=0$, we obtain the
full projectivized cotangent bundle. Of course, in any case $\tilde M$
carries a canonical contact structure and the analog of the Fefferman
construction refines this to a Lie contact structure. This generalizes
and explains the results of \cite{Sato-Yamaguchi}.

\noindent
(3) Consider the inclusion of $G:=Sp(2n,\Bbb R)$ into $\tilde
G:=SL(2n,\Bbb R)$ by the standard representation. Parabolic subgroups
of $G$ correspond to isotropic flags in the symplectic vector space
$\Bbb R^{2n}$, which parabolic subgroups in $\tilde G$ correspond to
arbitrary flags. Hence there is only one choice for a parabolic
subgroup $\tilde P\subset\tilde G$ such that the $G$--orbit of
$e\tilde P$ in $\tilde G/\tilde P$ is open. Namely, one has to use the
stabilizer of a line, since for lines being isotropic is a vacuous 
condition. In this case $P:=G\cap\tilde P$ is itself parabolic in
$G$. 

Hence we conclude that the analog of the Fefferman construction this
time starts from a geometry of type $(G,P)$ on $M$ and produces the
underlying structure of a geometry of type $(\tilde G,\tilde P)$ on
the same space $M$. Geometries of type $(G,P)$ are a contact analog of
projective structures, and our construction extends such a structure
to a classical projective structure. This has been directly obtained
in \cite{Fox}, where more details about such structures can be found. 

\noindent
(4) To finish, we discuss an exotic example which however has a long
history. Let $G$ be the split real form of the exceptional Lie group
$G_2$. It is well known that $G_2$ has a 7 dimensional representation,
and for the split form there is an invariant inner product of
signature $(3,4)$ on this representation. Hence this gives rise to an
inclusion of $G$ into $\tilde G:=SO(3,4)$. The stabilizer $P\subset G$
of a line through a highest weight vector in this representation is
one of the two maximal parabolic subgroups of $G$. This line is easily
seen to be isotropic, so as in (3) we obtain $P=G\cap\tilde P$, where
$\tilde P$ is the stabilizer of the highest weight line in $\tilde
G$. 

Geometries of type $(G,P)$ are exactly the generic rank two
distributions in dimension five which are studied in Cartan's famous
``five variables paper'' \cite{Cartan:five}. Given such a distribution
on $M$, the analog of the Fefferman construction produces a canonical
conformal class of split signature $(2,3)$ on $M$. Such a canonical
conformal class was recently discovered by P.~Nurowski using Cartan's
method of equivalence, see \cite{Nurowski}. Since in Nurowski's
construction one obtains the same normal Cartan connection for both
geometries, it is very likely that the structure described here
coincides with his.

\end{document}